\input epsf
\magnification=\magstep1
\baselineskip=15pt
\mathsurround=1pt

\font\Bbb=msbm10

\def\Z{\hbox{\Bbb Z}}
\def\Q{\hbox{\Bbb Q}}
\def\R{\hbox{\Bbb R}}
\def\C{\hbox{\Bbb C}}
\def\O{{\cal O}}
\def\H{{\cal H}}
\def\I{{\cal I}}

\centerline{\bf Bianchi Orbifolds of Small Discriminant}
\medskip 
\centerline{A. Hatcher}
\bigskip

Let $ \O_D $ be the ring of integers in the imaginary quadratic field $ \Q(\sqrt D) $ of discriminant $ D < 0 $. Then $ PGL_2(\O_D) $ is a discrete subgroup of the isometry group $ PSL_2(\C) $ ($=PGL_2(\C)$) of hyperbolic $3$-space $ \H^3 $. The quotient space $ \H^3/PGL_2(\O_D) = X_D $ is topologically a noncompact $3$-manifold whose cusps (ends) are of the form $ S^2 \times [0,\infty) $. The number of cusps of $ X_D $ is known to be $ h_D $, the class number of $ \O_D $. So $ X_D $ is a closed manifold $ \widehat X_D $ with $ h_D $ points removed.

For small $ D $, including the $ 31 $ discriminants in the range $ D > -100 $, R. Riley has done computer calculations of the Ford fundamental domain $ F_D $ for the action of $ PGL_2(\O_D) $ on $ \H^3 $. (See [5] for an account of the techniques; for about half of these $ D $'s, Bianchi had calculated the fundamental domains --- by hand, presumably --- almost a century ago, in [2],[3].) Riley's computer output includes how the faces of $ F_D $ are identified by elements of $ PGL_2(\O_D) $. So it becomes a pleasant exercise in geometric visualization to try to recognize the manifold $ \widehat X_D $. The results of carrying out this exercise for $ D > -100 $ are listed in Table I. ($ P^3 $ denotes real projective $ 3 $-space, $ \sharp $ is connected sum.)

\bigskip
\noindent\hskip18pt\vbox{\offinterlineskip
\halign{&\ \hfil#\hfil\ \cr
$D$\strut&\vrule&$-3$&$-4$&$-7$&$-8$&$-11$&$-15$&$-19$&$-20$&$-23$ &$-24$&$-31$&$-35$&$-39$ \cr
\noalign{\hrule}\cr
$\widehat X_D$\strut&\vrule height11pt &$S^3$&$S^3$&$S^3$ &$S^3$&$S^3$ &$S^3$&$S^3$&$S^3$
&$S^3$&$S^3$&$S^3$&$S^3$&$S^3$\cr }}
\medskip
\noindent\hskip50pt\vbox{\offinterlineskip
\halign{&\ \hfil#\hfil\ \cr
$-40$\strut&$-43$&$-47$&$-51$&$-52$&$-55$&$-56$&$-59$&$-67$&$-68$ &$-71$&$-79$\cr
\noalign{\hrule\vskip2pt}\cr
$P^3$\strut&$P^3$&$S^3$&$S^3$&$P^3$&$P^3$&$S^3$&$S^3$
&$P^3\sharp P^3$&$S^3$&$S^3$&$P^3$\cr }}
\medskip
\noindent\hskip50pt\vbox{\offinterlineskip
\halign{&\ \hfil#\hfil\ \cr
$-83$\strut&$-84$&$-87$&$-88$&$-91$&$-95$ \cr
\noalign{\hrule\vskip2pt}\cr
$P^3$\strut&$S^1\times S^2$&$S^1\times S^2$&$P^3\sharp P^3\sharp P^3$&$P^3\sharp P^3$&$P^3$\cr }}
\medskip
\centerline{Table I}
\medskip

There are exactly $ 19 $ $ D $'s in this range for which $ \widehat X_D $ is the $ 3 $-sphere. Since $ \pi_1X_D $ ($= \pi_1\widehat X_D $) is $ PGL_2(\O_D)/torsion $, the question of when $ \widehat X_D $ is $ S^3 $ is equivalent (assuming no $ \widehat X_D $'s are counterexamples to the Poincar\'e Conjecture) to when $ PGL_2(\O_D) $ is generated by torsion. In the $ 19 $ cases when $ \widehat X_D = S^3 $, we have determined in addition the orbifold structure on $ X_D $, namely, the embedded graph in $ X_D $ consisting of images under the quotient map $ \H^3 \to X_D $ of axes of rotations of torsion elements of $ PGL_2(\O_D) $, each edge of this graph being labelled by the order of the corresponding torsion element. These orbifold structures are shown in the figures on the next page. (Labels ``2" on edges are omitted. The small circles denote the cusp spheres. The sphere $ S^3 = \widehat X_D $ is regarded as $ 3 $-space compactified by a point at infinity.)

\eject

\centerline{\epsfbox{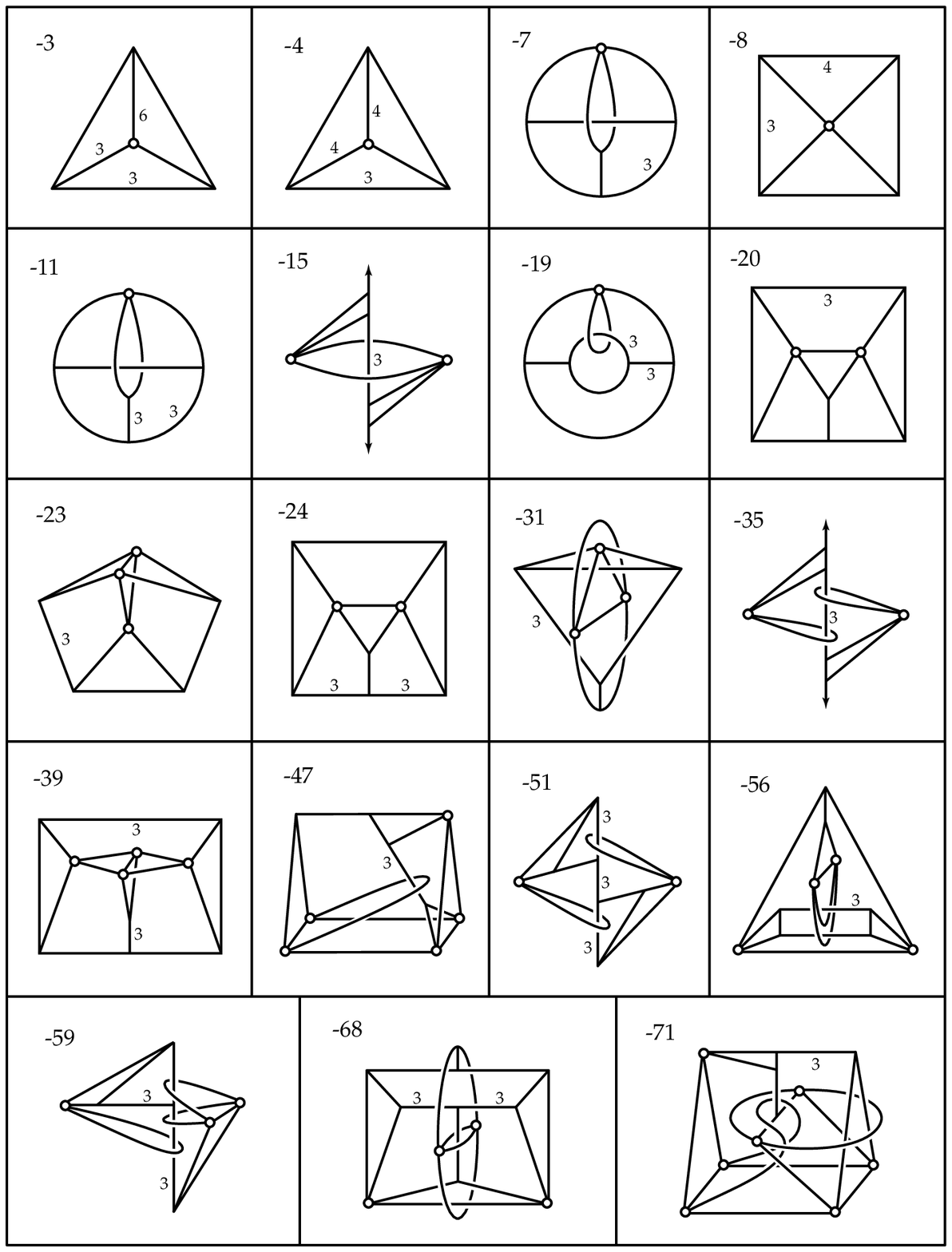}}
\eject

Table II below shows the orbifold structure on $ X_D $ in a few cases when $ \widehat X_D \ne S^3 $. In the top row are the first four cases when $ \widehat X_D = P^3 $. Here we view $ P^3 $ as a ball with antipodal points of its boundary sphere (indicated by the dashed-line circles) identified. The lower part of the Table represents the case $ D = -84 $, when $ \widehat X_D = S^1 \times S^2 $. The periodic extension of the graph shown, modulo its translation symmetries, gives the ``singular" locus of the orbifold structure on $ X_{-84} $, a graph lying on the torus $ S^1 \times S^1 $ which decomposes $ S^1 \times S^2 $ into two $ S^1 \times D^2 $'s. The vertical direction in the figure represents the meridian circles $ \{ x\} \times \partial D^2 $ on these $ S^1 \times D^2 $'s. 

\medskip
\centerline{\epsfbox{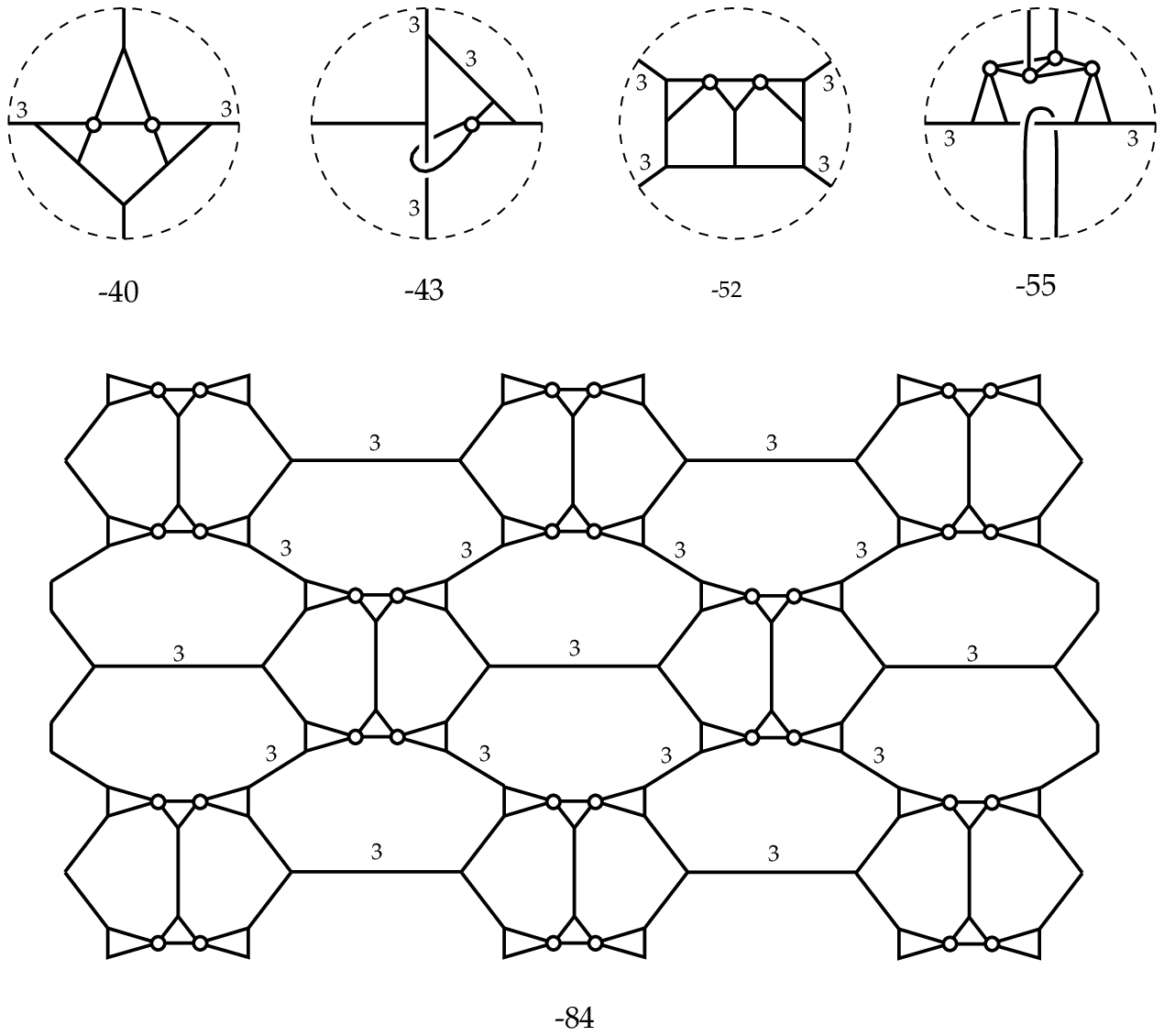}}
\centerline{Table II}
\medskip

For the index two subgroup $ PSL_2(\O_D) $ of $ PGL_2(\O_D) $,  the quotient space $ Y_D
=
\H^3/PSL_2(\O_D) $ is a $ 2 $-sheeted branched cover of $ X_D $, branched in such a way
that the cusp spheres of $ X_D $ become cusp tori of  $ Y_D $ (except for $ D = -3,-4
$, when they remain spheres). It turns out that in the $ 19 $ cases when $ \widehat X_D
= S^3 $, this branching condition at cusps uniquely determines the branched covering $
Y_D \to X_D $, and one can very easily by inspection determine the topological type of $
Y_D $. This is given in Table III, in which $ Y_D $ is strictly speaking the interior of
the compact manifold listed. (Notations: $ B^3 = $ 3-ball, $ D^2 = $ 2-disk, $ T^2 = $
torus, $ I = [0,1] $.)

\bigskip
\noindent\hskip15pt\hbox{\vbox{\offinterlineskip
\vbox{\halign{&\ #\hfil\ \cr
\ $D$\strut&\vrule&$Y_D$\cr
\noalign{\hrule}\cr
$-3$\strut&\vrule height11pt&$B^3$\cr
$-4$\strut&\vrule&$B^3$\cr
$-7$\strut&\vrule&$S^1\times D^2$\cr
$-8$\strut&\vrule&$S^1\times D^2$\cr
$-11$\strut&\vrule&$S^1\times D^2$\cr
$-15$\strut&\vrule&$S^1\times D^2\ \sharp\ S^1\times D^2$\cr
$-19$\strut&\vrule&$S^1\times D^2$\cr
$-20$\strut&\vrule&$S^1\times D^2\ \sharp\ S^1\times D^2$\cr
$-23$\strut&\vrule&$S^1\times D^2\ \sharp\ T^2\times I$\cr
$-24$\strut&\vrule&$S^1\times D^2\ \sharp\ S^1\times D^2$\cr  }}}

\hskip20pt\hbox{\vbox{\offinterlineskip
\halign{&\ #\hfil\ \cr
\ $D$\strut&\vrule&$Y_D$\cr
\noalign{\hrule}\cr
$-31$\strut&\vrule height11pt&$S^1\times D^2\ \sharp\ T^2\times I$\cr
$-35$\strut&\vrule&$S^1\times D^2\ \sharp\ S^1\times D^2 \ \sharp\ S^1\times S^2$\cr
$-39$\strut&\vrule&$S^1\times D^2\ \sharp\ S^1\times D^2 \ \sharp\ T^2\times I$\cr
$-47$\strut&\vrule&$S^1\times D^2\ \sharp\ T^2\times I \ \sharp\ T^2\times I$\cr
$-51$\strut&\vrule&$S^1\times D^2\ \sharp\ S^1\times D^2 \ \sharp\ S^1\times S^2$\cr
$-56$\strut&\vrule&$S^1\times D^2\ \sharp\ S^1\times D^2 \ \sharp\ T^2\times I \ \sharp\ S^1 \times S^2$\cr
$-59$\strut&\vrule&$S^1\times D^2 \ \sharp\ T^2\times I \ \sharp\ S^1 \times S^2$\cr
$-68$\strut&\vrule&$S^1\times D^2\ \sharp\ S^1\times D^2 \ \sharp\ T^2\times I \ \sharp\ S^1 \times S^2$\cr
$-71$\strut&\vrule&$S^1\times D^2\ \sharp\ T^2\times I \ \sharp\ T^2\times I \ \sharp\
T^2\times I$\cr
\strut& &\cr  }}}}

\medskip
\centerline{Table III}
\medskip

In the $ 14 $ cases that $ Y_D $ does not contain a  connected summand $ S^1 \times S^2
$, the restriction map $ H^1(Y_D;\Z) \to H^1(\partial Y_D;\Z) $ is injective; in other
words, ``$ Y_D $ has no cuspidal cohomology." It is known (see [9],[4],[1],[6],[8]) that
these are the only cases when this happens, for arbitrary $ D < 0 $.

Perhaps the first thing one notices about the pictures of the orbifolds $ X_D $  is the
symmetries. In each case there is a reflectional symmetry through a plane parallel to
the plane of the page. This reflection presumably corresponds to the $ \Z_2 $ extension
of $ PGL_2(\O_D) $ obtained by adjoining complex conjugation, the Galois automorphism of
$ \O_D $. When $ D $ has more than one distinct prime divisor, there is also a $
180^\circ $ rotational symmetry evident in the pictures. (This symmetry does not appear
in the Ford domains, however.) Such a symmetry is predicted by general theory: Bianchi
[3] already described a group $ G_D \subset PGL_2(\C) $ containing $ PGL_2(\O_D) $ as a
normal subgroup of finite index, with quotient $ G_D/PGL_2(\O_D) \approx \I_D/2\I_D $,
the mod $2 $ ideal class group (or genera group) of $ \O_D $. According to Gauss, this
quotient has $ \Z_2 $-rank equal to one less than the number of distinct prime divisors
of $ D $. 

Even when $ \widehat X_D $ is not $ S^3 $,  $ \widehat X_D $ may have $ S^3 $ as the
quotient space corresponding to a finite extension of $ PGL_2(\O_D) $, such as the group
$ G_D $ above. This happens in the cases $ D = -40,-52,-55 $ in Table II, when $
\widehat X_D = P^3 $. It also happens for $ D = -84 $, when $ \widehat  X_D = S^1
\times S^2 $ has a $ \Z_2 \times \Z_2 $ ($\approx \I_{-84}/2\I_{-84} $) quotient which
is $ S^3 $, the $ \Z_2 \times \Z_2 $ action on $ X_{-84} $ restricting to the full
symmetry group of the singular locus on the torus shown in Table II.

It appears that except for $ D = -3,-4 $, the remaining  $ 17 $ orbifolds $ X_D $ with $
\widehat X_D = S^3 $ are Haken orbifolds [7]. That is, by repeatedly splitting open
along incompressible $ 2 $-dimenional suborbifolds, $ X_D $ can be reduced to a disjoint
union of finitely many orbifolds of the form $ \R^3/\Gamma $ for $ \Gamma $ a finite
subgroup of $ SO(3) $ (acting on $ \R^3 $ as isometries). These splitting surfaces are
all separating, so such a hierarchy for $ X_D $ yields a way of building up $
PGL_2(\O_D) $ from finite subgroups of $ SO(3) $ by iterated free product with
amalgamation constructions. These hierarchies are in general far from unique. As a very
simple example, the orbifold $ X_{-8} $ as drawn can be split successively along
horizontal and vertical planes through the cusp, in either order, yielding the two
structures
$$
\bigl(O(24) *_{C(4)} D(8)\bigr)\ *_{C(3)*C(2)}\ \bigl(D(6)*_{C(2)}D(4)\bigr)
$$
and 
$$
\bigl(O(24) *_{C(3)} D(6)\bigr)\ *_{C(4)*C(2)}\ \bigl(D(8)*_{C(2)}D(4)\bigr)
$$
where $ O(24) $ is the octahedral group, $ D(2n) $ is the dihedral group of order $ 2n $, and $C(n) $ is the cyclic group of order $ n $. In more subtle examples, not even the collection of finite subgroups of $ SO(3) $ which start the iterated amalgamated free product construction is unique, though of course the noncyclic subgroups among these are unique, corresponding to the vertices in the singular locus of the orbifold structure.

In all cases except $ D = -3,-4 $, there is a splitting
$$
PGL_2(\O_D) \approx PGL_2(\Z) *_A(?)
$$
amalgamated over $ A = PSL_2(\Z) $, arising as follows.  In the upper half-space model
of $ \H^3 $, bounded below by the plane $ \C $, there lies $ \H^2 $, the half-plane
above $ \R $. The orbifold $ \H^2/PGL_2(\Z) $ is a triangle with one vertex at the cusp
at $ \infty $. This triangle is embedded in $ X_D $, and the boundary of a small regular
neighborhood of this triangle is the surface corresponding to $ PSL_2(\Z) $ in the
splitting above. This surface can be take to be totally geodesic in $ X_D $. It should
be of interest to find other totally geodesic incompressible surfaces in $ X_D $, since
these are more likely to be defined arithmetically. For example, as Riley has pointed
out, the cuspidal classes in $ H^1(Y_D;\Z) \approx H^2(Y_D,\partial Y_D;\Z) $ found in
[9],[4],[1] are represented by totally geodesic surfaces formed by the intersections of
the Ford domain with certain planes parallel to $ \H^2 \subset \H^3 $. These
non-separating surfaces in $ Y_D $ pass down to non-separating (totally geodesic)
surfaces in $ X_D $, which are often non-orientable. Since non-separating surfaces do not
exist in
$ S^3 $, it follows from [9],[4],[1] that the only values of $ D < -100 $ for which $
\widehat X_D $ could be $ S^3 $ are $ -119,-164, -191, -311,-356, -404,-479$, and
$-776 $. Riley's computer calculations eliminate $ -164 $ from this list.

\bigskip\noindent
{\bf References}
\medskip
\item{[1]} M.Baker, Ramified primes and the homology of the Bianchi groups, I.H.E.S. preprint (1982).

\item{[2]} L. Bianchi, Sui gruppi di sostitutioni lineari con coefficienti  a corpi
quadratici imaginarii, {\it Math. Annalen} 40 (1892), 332-412. [This article and the
following one are reprinted in Bianchi's collected workes, {\it Opere}, vol.1, Editioni
Cremonese, Rome, 1952.]

\item{[3]} L. Bianchi, Sui gruppi di sostitutioni lineari, {\it Math. Annalen} 42 (1892),
30-57.

\item{[4]} F. Grunewald and J. Schwermer, Arithmetic quotients of hyperbolic 3-space, 
cusp forms and link complements, {\it Duke Math. J.} 48 (1981), 351-358.

\item{[5]} R. Riley, Application of a computer implementation of Poincar\'e's  theorem
on fundamental polyhedra, {\it Math. of Computation} 40 (1983), 607-632.

\item{[6]} J. Rohlfs, On the cuspidal cohomology of the Bianchi modular groups, {\it Math. Z.} 188 (1985), 253-269.

\item{[7]} W. Thurston, Geometry and topology of 3-manifolds, xeroxed notes.

\item{[8]} K. Vogtmann, Rational homology of Bianchi groups, {\it Math. Annalen} 272 (1985), 399-419.

\item{[9]} R. Zimmert, Zur $SL_2$ der ganzen  Zahlen eines imagin\"ar quadratischen
Zahlk\"orpers, {\it Invent. math.} 19 (1973), 73-82.

\medskip\noindent
December 1983

\noindent
Cornell University

\end